\documentclass[12pt]{amsart}

\newcommand{\vsp}{\vspace{5mm}}
\newcommand{\bt}{\begin{Theorem}}
\newcommand{\et}{\end{Theorem}}
\newcommand{\bi}{\beign{itemize}}
\newcommand{\ei}{\end{itemize}}
\newcommand{\bea}{\begin{eqnarray}}
\newcommand{\eea}{\end{eqnarray}}
\newtheorem{Theorem}{\sc Theorem}
\newtheorem{Lemma}[Theorem]{\sc Lemma}
\newtheorem{Proposition}[Theorem]{\sc Proposition}
\newtheorem{Corollary}[Theorem]{\sc Corollary}
\newtheorem{Definition}[Theorem]{\sc Definition}
\newtheorem{Example}[Theorem]{\sc Example}
\newtheorem{Remark}[Theorem]{\sc Remark}

\newcommand{\be}{\begin{equation}}
\newcommand{\ee}{\end{equation}}

\def\qed{\hfill$\diamondsuit$}
\def\bi{\bibitem}

\def\CM{{\mathcal {M}}}
\def\CB{{\mathcal {B}}}

\def\CL{{\mathcal {L}}}
\def\CH{{\mathcal {H}}}
\def\CK{{\mathcal {K}}}

\def\CS{{\mathcal {S}}}

\def\CD{{\mathcal {D}}}
\def\CO{{\mathcal {O}}}
\def\CZ{{\rm\ke}rn.26em
\newcommand\la{{\langle}}
\newcommand\ra{{\rangle}}
\newcommand\lar{\leftarrow}
\newcommand\Lar{\Leftarrow}
\newcommand\rar{\rightarrow}
\newcommand\Rar{\Rightarrow}

\vrule width.02em height.5ex depth0ex \kern.04em \vrule width
.02em height1.47ex depth-1ex \kern-.34em Z}

\def\C{{\rm \kern.24em
 \vrule width.02em
    height1.4ex depth-.05ex
 \kern-.26em C}}

\def\ra{{\rightarrow}}

\def\ei{{\bf e_i}}

\def\us{{\underline{s}}}

\def\uT{{\underline{T}}}
\def\uR{{\underline{R}}}
\def\uV{{\underline{V}}}
\def\uW{{\underline{W}}}
\def\uS{{\underline{S}}}

\def\uZ{{\underline{Z}}}

\def\circleds{{\bigcirc \!\!\!\!s}}
\def\N{{\rm I\kern-.23em N}}
\def\B{{\rm I\kern-.25em B}}
\def\D{{\rm I\kern-.25em D}}
\def\E{{\rm I\kern-.25em E}}
\def\F{{\rm I\kern-.25em F}}

\def\I{{\rm I\kern-.25em I}}

\def\M{{\rm I\kern-.23em M}}
\def\P{{\rm I\kern-.25em P}}
\def\A{{\rm \kern.22em
 \vrule width.02em
    height0.5ex depth 0ex
 \kern-.24em A}}
\def\G{{\rm \kern.24em
 \vrule width.02em
    height1.4ex depth-.05ex
 \kern-.26em G}}
\def\J{{\rm \kern.19em
 \vrule width.02em
    height1.47ex depth 0ex
 \kern-.21em J}}
\def\O{{\rm \kern.24em
 \vrule width.02em
    height1.4ex depth-0.5ex
 \kern-.26em O}}
\def\Q{{\rm \kern.24em
 \vrule width.02em
    height1.4ex depth-.05ex
 \kern-.26em Q}}
\def\S{{\rm \kern.18em
 \vrule width.02em
    height1.4ex depth-.9ex
  \kern.12em
  \vrule width.02em
     height0.7ex depth 0ex
  \kern-.34em S}}
\def\T{{\rm \kern.45em
 \vrule width.02em
    height1.47ex depth 0ex
 \kern-.47em T}}
\def\U{{\rm \kern.30em
 \vrule width.02em
    height1.47ex depth-.05ex
 \kern-.32em U}}
\def\V{{\rm \kern.27em
 \vrule width.02em
    height1.47ex depth-.8ex
 \kern-.29em V}}
\def\W{{\rm \kern.25em
 \vrule width.02em
    height1.47ex depth-0.9ex
 \kern.34em
 \vrule width.02em
    height1.47ex depth-.9ex
  \kern-.63em W}}
\def\X{{\rm \kern.30em
 \vrule width.02em
    height1.4ex depth-1ex
  \kern.12em
  \vrule width.02em
     height0.4ex depth 0ex
  \kern-.46em X}}
\def\Y{{\rm \kern.25em
 \vrule width.02em
    height1.0ex depth 0ex
 \kern-.27em Y}}
\def\Z{{\rm \kern.26em
 \vrule width.02em
    height0.5ex depth 0ex
  \kern.04em
  \vrule width.02em
     height1.47ex depth-1ex
  \kern-.34em Z}}

\begin{document}

\begin{center} {\bf {\Large Standard noncommuting and commuting dilations
of commuting tuples}}
\end{center}

\vsp
\vsp
\vsp

\begin{center}

{\sc B. V. Rajarama Bhat,  Tirthankar Bhattacharyya, and Santanu Dey }
\vsp

{\bf April 5, 2002}
\end{center}

\vsp \vsp \vsp \vsp

\begin{center}
{\underline {\bf Abstract}}
\end{center}

\vsp
We introduce a notion called  `maximal commuting piece' for tuples of Hilbert
space operators. Given a commuting tuple of operators forming a row
contraction there are two commonly  used dilations in multivariable
operator theory. Firstly there is the minimal isometric
dilation consisting of isometries with orthogonal ranges and hence it
is a noncommuting tuple. There is also a  commuting dilation related
with a standard commuting tuple on Boson Fock space. We show
that this commuting dilation is the maximal commuting piece of the
minimal isometric dilation.  We use this result to classify all
representations of Cuntz algebra $\CO _n$ coming from dilations of commuting tuples.

\vsp
\vfill
----------------------------------------------------------------------



\newpage

\begin{section}{Introduction}
\setcounter{equation}{0}
It is a well-known result due to Sz.\@ Nagy that every contraction
on a Hilbert space dilates to an isometry. There is a very natural
generalization of this result to a class of operator tuples defined as
follows.
\begin{Definition}{\em
 A {\em contractive
$n$-tuple,\/} or a {\em row contraction\/}
 is a $n$-tuple $\underline{T}=(T_1,\ldots ,T_n)$ of bounded
operators on a Hilbert space ${\mathcal H}$ such that $T_1T_1^*+
\cdots  +T_nT_n^*\leq I$. }
\end{Definition}

 Such tuples are known as row contractions as the condition is equivalent to
having the operator $(T_1, \ldots T_n)$ from
$\CH \oplus \cdots \oplus \CH$ ($n$-times)
to $\CH $,  as a contraction. It is possible to dilate contractive tuples to tuples
of isometries with orthogonal ranges. Moreover, such a dilation is
unique up to unitary equivalence, under a natural minimality condition
just as in the one variable case. This dilation which we call as {\em
the minimal isometric dilation\/ } or as {\em standard noncommuting
dilation\/ } has been explored by many authors.
Some ideas along this direction can already be seen in the early
paper \cite{Da} of Davis. In more concrete form this dilation can be seen
in the papers of Bunce \cite{Bu} and Frazho \cite{Fr1, Fr2}. A real extensive study of
this notion has been carried out by Popescu in a series of papers (\cite{Po1}-\cite{Po5}, \cite{AP1, AP2})
 and he has neat generalizations of many results from one variable
situation.

Now suppose the tuple under consideration is a commuting tuple in the sense that
$T_iT_j=T_jT_i$ for all $1\leq i, j\leq n.$ Then it is natural to wish for
a dilation consisting of mutually commuting isometries. Unfortunately,
such a dilation does not exist in general for $n\geq 3$ \cite{Pa}.
However, there is a dilation of commuting contractive tuples studied
recently by Arveson \cite{Ar2}. In fact, this particular dilation was  looked at by
Drury \cite{Dr} in his study of von Neumann inequality for tuples and similar ideas
have been explored by Agler \cite{Ag}, Athavale \cite{At2} and
others  for different classes of operators
using various reproducing kernels. We call this
dilation of Drury and Arveson  as {\em standard commuting dilation\/ } of contractive commuting tuples.
This dilation consists of a commuting tuple, but the constituent operators
are not isometries. Then a natural question arises as to in what sense this
dilation is canonical. Here in Section 3 as our main
result, we show that the standard commuting dilation is the
 `maximal commuting tuple contained' in
the minimal isometric dilation.  To begin with we make these words inside inverted
commas precise
by defining what we call as `maximal commuting piece' of a tuple of operators.
We can explore how the standard commuting dilation of the maximal commuting piece sits
inside the minimal isometric dilation of the original tuple, and whether it
is the maximal commuting piece of the minimal isometric dilation etc.
 We have been able to carry out this study for purely contractive tuples in Section 2.

Any tuple $(W_1, \ldots W_n)$ of isometries with orthogonal ranges satisfying
$\sum W_iW_i^*=I$ gives us a concrete representation of the familiar Cuntz
algebra \cite{Cu}. Recently, there has been a lot of effort to study such representations
in connection with wavelet theory, see for instance the
papers  \cite{BJ1, BJ2} of Bratteli and Jorgensen. If we start with
a contractive tuple $(T_1, \ldots , T_n)$ satisfying $\sum T_iT_i^*=I$ and
consider the minimal isometric dilation we actually have a representation
of the Cuntz algebra. Very interesting results on classification of
these representation up to unitary
equivalence in terms of invariants determined by $(T_1, \ldots T_n)$ have
been obtained by Davidson, Kribs, and Shpigel \cite{DKS}, where the
operators $T_i$ act
 on a finite dimensional space. It is a natural question as
to what are the representations of Cuntz algebra one can get by dilating
contractive tuples which are also commuting. Surprisingly they are very
few and are all determined by the GNS representations of the so-called Cuntz states. This result we obtain
in Section 4, as an application of the  main result. Unlike the work of
Davidson and Kribs, we need not restrict  the operators $T_i$ to be
acting on finite dimensional spaces.

All the Hilbert spaces we consider will be complex and separable.  For
a subspace $\CH $ of a Hilbert space $P_{\CH }$ will denote the orthogonal
projection onto $\CH .$ For fixed
$n\geq 2,$ we need
two standard $n$-tuples of operators,  denoted by $\uV $ and $\uS$
acting on Fock spaces. For any Hilbert space $\CK $, we have the full Fock
space over $\CK$ denoted by $\Gamma (\CK)$ and the Boson (or symmetric) Fock
space over $\CK$ denoted by $\Gamma _s(\CK)$ as,
$$\Gamma (\CK)=\mathbb{C}\oplus \CK \oplus \CK^{\otimes ^2}\oplus \cdots
\oplus \CK ^{\otimes ^m}\oplus \cdots ,
$$
$$\Gamma _s(\CK)=\mathbb{C}\oplus \CK\oplus \CK ^{\circleds ^2}\oplus \cdots
\oplus \CK ^{\circleds ^m}\oplus \cdots,$$
where $\CK ^{\circleds ^m}$ denotes $m$-fold symmetric tensor product.
We will consider the Boson Fock space as a subspace of the full Fock
space in the natural way. We denote the vacuum vector $1\oplus 0\oplus \cdots $
(in either of
Fock spaces) by $\omega .$
Let ${\mathbb{C}}^n$ be the $n$-dimensional complex
Euclidian space with usual inner product and $\Gamma({\mathbb{C}}^n)$ be the full Fock
space over ${\mathbb{C}}^n$. Let $\{e_1, \ldots , e_n\}$ be the standard ortho-normal
basis of ${\mathbb{C}}^n$. Then the (left) creation operators $V_i$ on
$\Gamma({\mathbb{C}^n})$ are defined by
\[ V_i x = e_i \otimes x \]  where $1 \leq i \leq n$ and $x
\in
\Gamma({\mathbb{C}}^n)$ (Of course, here $e_i\otimes \omega$ is interpreted
as $e_i$).
It is obvious  that the tuple $\underline{V}= (V_1, \ldots , V_n)$  consists
of isometries with orthogonal ranges and it is contractive, in fact
$\sum V_iV_i^*=I-E_0$, where $E_0$ is the projection on to the vacuum space.
Let $\uS =(S_1, \ldots , S_n)$ be the tuple of operators on $\Gamma _s(\mathbb{C}^n)$
where, $S_i$ is the compression of $V_i$ to $\Gamma _s(\mathbb{C}^n)$:
$$S_i=P_{\Gamma _s(\mathbb{C}^n)}V_i|_{\Gamma _s(\mathbb{C}^n)}.$$
Clearly each  $V_i^*$ leaves $\Gamma _s(\mathbb{C}^n)$ invariant.
Therefore, $S_i^*x=V_i^*x$, for $x\in \Gamma _s(\mathbb{C}^n)$. Then it is
easy to see that $(S_1, \ldots S_n)$ is also a contractive tuple satisfying
$\sum S_iS_i^*=I^s-E^s_0$ (where now $I^s, E^s_0$ are  identity, projection onto vacuum
space in $\Gamma _s(\mathbb{C}^n)$. Moreover a simple computation
shows that $S_iS_j=S_jS_i$ for all $1\leq i, j\leq n.$

For  operator tuples $(T_1, \ldots , T_n)$ quite often we need to consider
 products of the form $T_{\alpha _1}T_{\alpha _2}\cdots$ $ T_{\alpha _m}$,
where each $\alpha _k\in \{ 1, 2, \ldots , n\}$. So it is convenient to
have a notation for such products. Let $\Lambda $ denote the set
$\{ 1, 2, \ldots , n\}$ and $\Lambda ^m$ denote the $m$-fold cartesian
product of $\Lambda $ for $m\geq 1.$ Given $\alpha =(\alpha _1, \ldots , \alpha _m)$
in $\Lambda ^m$, $\uT ^{\alpha }$ will mean the operator
$T_{\alpha _1}T_{\alpha _2}\cdots T_{\alpha _m}$.
Let $\tilde{\Lambda}$ denote $\cup_{n=0}^{\infty} \Lambda^n$, where
$\Lambda ^0$ is just the set $\{ 0\}$ by convention and by $\uT ^0$ we
would mean the identity operator of the  Hilbert space where the operators
$T_i$
are acting. In a similar fashion for $\alpha \in {\tilde \Lambda}$,
$e^{\alpha }$ will denote the vector $e_{\alpha _1}\otimes e_{\alpha _2}\otimes
\cdots \otimes e_{\alpha _m}$ in  the full Fock space $\Gamma (\mathbb{C}^n)$
and $e^0$ is the vacuum $\omega $.

\end{section}

\begin{section}{Maximal Commuting Piece and Dilation}
\setcounter{equation}{0}

\begin{Definition}{\em
Let ${\mathcal H}, {\mathcal L}$ be two Hilbert spaces such that
${\mathcal H}$ is a closed subspace of ${\mathcal L}$. Suppose
$\uT, \uR$ are $n$-tuples of bounded operators on ${\mathcal H}$,
${\mathcal L}$ respectively. Then $\uR $ is called a {\em dilation\/} of
 $\uT $ if
$$R_i^*u=T_i^*u $$
 for all $u\in {\mathcal H}, 1\leq i\leq n.$
In such a case $\uT$ is called a {\em piece\/} of $\uR .$
If further $\uT$ is a commuting tuple ( i.\@ e., $T_iT_j=T_jT_i,$ for
all $i,j$), then it is called a commuting piece of $\uR .$ A dilation
$\uR $ of $\uT$ is said to be a {\em minimal dilation\/ } if
${\overline {\mbox span}}\{ \uR ^{\alpha }h: \alpha \in {\tilde \Lambda},
h\in \CH \}=\CL .$}
\end{Definition}

In this Definition we note that if $\uR$ is a dilation of $\uT$,
then ${\mathcal H}$ is a co-invariant subspace of $\uR$, that is,
it is left invariant by all $R_i^*$. It is standard (see \cite{Ha})
to call $(R_1^*, \ldots ,
R_n^*)$ as an {\em extension\/} of $(T_1^*, \ldots , T_n^*)$ and
$(T_1^*, \ldots , T_n^*)$ as a {\em part\/} of
$(R_1^*, \ldots , R_n^*)$. In such a situation it is easy to see
that  for any $\alpha, \beta \in {\tilde \Lambda },$
$\uT ^{\alpha }(\uT ^{\beta })^*$ is the compression of
$\uR ^{\alpha }(\uR ^{\beta })^*$ to ${\mathcal H}$, that is,
\begin{equation}
\uT ^{\alpha }(\uT ^{\beta })^*= P_{{\mathcal H}}\uR ^{\alpha }(\uR ^{\beta
})^*|_{{\mathcal H}}.
\end{equation}
We may extend this relation to any polynomials $p, q$ in
$n$-noncommuting variables to have
$$p(\uT)(q(\uT))^*=P_{{\mathcal H}}p(\uR )(q(\uR))^*|_{{\mathcal
H}}.$$
Usually it is property (2.1) is all that one demands of a dilation. But
we have imposed a condition of co-invariance in Definition 2, as it
is very convenient to have it this way for our purposes.

Now we look at commuting pieces of tuples.
Let $\uR $ be a $n$-tuple of bounded operators on a Hilbert space
$\CL$. Consider
 $${\mathcal C}(\uR )=\{ \CM : \CM ~~\mbox {is a co-invariant subspace
 for each }~~R_i, ~~R_i^*R_j^*h=R_j^*R_i^*h, \forall h\in \CM, \forall i,j\}.$$
So ${\mathcal C}(\uR)$ consists  of all co-invariant
subspaces of a $n$-tuple   of operators $\uR $
such that the compressions form a commuting tuple. It is a complete
lattice, in the sense that arbitrary intersections and closures
of arbitrary unions of such spaces are again in this collection.
Therefore it has a maximal element. We denote it by $\CL ^c(\uR)$
(or by $\CL^c $ when the tuple under consideration is clear).

\begin{Definition}{\em

Suppose $\uR $ is a $n$-tuple of operators on a Hilbert space $\CL$.
Then the {\em
maximal commuting piece\/ } of $\uR$ is defined as
the commuting piece $\uR ^c=(R_1^c, \ldots , R_n^c)$ obtained by compressing $\uR$ to the maximal element
$\CL ^c(\uR)$
of ${\mathcal C}(\uR)$. The maximal commuting piece is said to be {\em trivial}
if the space $\CL ^c(\uR)$ is just the zero space.}
\end{Definition}

It is quite easy to get tuples with trivial commuting piece, as tuples with
no non-trivial co-invariant subspaces have this property.
Of course, our main interest lies in tuples with non-trivial commuting pieces.
The following result is quite useful in determining the maximal commuting
piece.
\begin{Proposition}
Let $\uR$ be a $n$-tuple of bounded operators on a Hilbert space
${\mathcal L}$. Let $\CK _{ij}=
\overline{span}\{ \uR ^{\alpha }(R_iR_j-R_jR_i)h: h\in \CL , \alpha
\in \tilde {\Lambda }\}$ for all $1\leq i, j\leq n$, and
$\CK=\overline{span}
\{\cup_{i, j =1}^n \CK _{ij} \}$.
Then $\CL ^c(\uR) = \CK ^{\perp }.$ In other words,
$\CL ^c(\uR)=\{ h\in \CL: (R_i^*R_j^*-R_j^*R_i^*)(\uR ^\alpha )^*h=0,
\forall 1\leq i, j\leq n, \alpha\in {\tilde \Lambda}\}.$
\end{Proposition}

\noindent {\sc Proof:} Firstly $\CK ^{\perp }$ is a co-invariant subspace
of $\uR$  is obvious as each $R_i$ leaves $\CK $
invariant. Now for $i,j \in \{1,2, \ldots , n\}
$, and $h_1 \in \CL ^c,  h_2\in \CL$,
 $$\langle (R_i^*R_j^*-R_j^*R_i^*)h_1, h_2\rangle
 =  \langle h_1, (R_jR_i-R_iR_j)h_2\rangle =0.$$
 So we get $(R_i^*R_j^*-R_j^*R_i^*)h_1=0$.
 Now if $\CM $ is an element of ${\mathcal C}(\uR)$,
  take $i, j\in \{1, \ldots , n\},
 \alpha \in \tilde {\Lambda }, h_1\in \CM , h\in \CL $. We
 have
 $$ \langle h_1, \uR ^{\alpha }(R_iR_j-R_jR_i)h\rangle
 =\langle (R_j^*R_i^*-R_i^*R_j^*)(\uR ^{\alpha })^*h_1, h\rangle
 =0$$
 as $(\uR ^{\alpha })^*h_1\in \CM .$ Hence $\CM $ is contained in
 $\CK ^{\perp } $ Now the last statement is easy to see.
 \qed

\begin{Corollary}
Suppose $\uR $, $\uT $ are $n$-tuples of operators on two Hilbert
spaces $\CL , \CM$. Then the maximal commuting piece of $(R_1\oplus T_1, \ldots ,
R_n\oplus T_n)$ acting on $\CL\oplus \CM$ is $(R_1^c\oplus T_1^c,
\ldots , R_n^c\oplus T_n^c)$ acting on $\CL ^c\oplus \CM ^c.$ The maximal commuting piece of
$(R_1\otimes I, \ldots , R_n\otimes I)$ acting on $\CL \otimes \CM$
is $(R_1^c\otimes I, \ldots , R_n^c\otimes I)$ acting on $\CL ^c\otimes \CM.$
\end{Corollary}

\noindent {\sc Proof:} Clear from Proposition 4.
\qed

\begin{Proposition}
Let  $\uV =(V_1, \ldots , V_n)$ and $\uS=(S_1, \ldots S_n)$ be
standard contractive tuples on full Fock space $\Gamma ({\mathbb{C}}^n)$
and Boson Fock space $\Gamma _s({\mathbb{C}}^n)$ respectively. Then the
maximal commuting piece of $\uV $ is $\uS$.
\end{Proposition}

\noindent {\sc Proof:} As we have already noted in the Introduction, $\uS$ is a commuting
piece of $\uV$. To show maximality we make use of Proposition 4.
Suppose $x \in \Gamma ({\mathbb{C}}^n)$ and $\langle x,
\uV ^\alpha(V_jV_j-V_jV_i)y\rangle =0$ for all $\alpha \in \tilde{\Lambda },
1\leq i,j\leq n$ and $y\in \Gamma ({\mathbb{C}}^n)$. We wish
to show that $x\in \Gamma _s({\mathbb{C}}^n)$. Suppose $x_m$ is
the $m$-particle component of $x$, that is, $x=\oplus _{m\geq 0}x_m$
with $x_m \in (\mathbb{C}^n)^{{\otimes}^m}$ for $m\geq 0$.
For $m\geq 2$ and any permutation $\sigma $ of $\{1,2, \ldots ,m\}$ we
need to show that the unitary $U_{\sigma }:(\mathbb{C}^n)^{{\otimes}^m}
\to (\mathbb{C}^n)^{{\otimes}^m}$, defined by
$$U_{\sigma }(u_1\otimes \cdots \otimes u_m)=
u_{\sigma ^{-1}(1)}\otimes \cdots \otimes u_{\sigma ^{-1}(m)},$$
leaves $x_m$ fixed. Since the group of permutations of $\{1, 2, \ldots ,m\}$
is generated by permutations $\{(1,2), \ldots , (m-1, m)\}$
it is enough to verify
$U_{\sigma }(x_m)=x_m$ for permutations
$\sigma $ of the form $(i, i+1)$. So fix $m$ and $i$ with $m\geq 2$ and
$1\leq i\leq (m-1).$  We have
 $$ \langle \oplus _px_p, \uV ^\alpha(V_kV_l-V_lV_k)y\rangle =0,$$
for every $y\in \Gamma (\mathbb{C}^n)$, $1\leq k,l\leq n.$ As $\alpha
$ is arbitrary, this means that
$$\langle x_m, z\otimes (e_k\otimes e_l-e_l\otimes e_k)\otimes w\rangle =0$$
for any $z\in (\mathbb{C}^{n})^{\otimes ^{(i-1)}}, w\in
 (\mathbb{C}^n)^{\otimes ^{(m-i-1)}}.$
This clearly implies $U_{\sigma }(x_m)=x_m$, for $\sigma =(i, i+1)$.
\qed

Now let us see how the maximal commuting piece behaves with respect
to the operation of taking dilations. Before considering specific
dilations we have the following general statement.
\begin{Proposition}
 Suppose
$\uT, \uR$ are $n$-tuples of bounded operators on ${\mathcal H}$,
${\mathcal L}$, with $\CH \subseteq \CL$, such that $\uR$ is
a dilation of $\uT$. Then $\CH ^c(\uT)=\CL ^c(\uR)\bigcap \CH$ and
$\uR ^c$ is a dilation of $\uT ^c$.
\end{Proposition}

\noindent {\sc Proof:} We have $R_i^*h=T_i^*h,$ for $h\in \CH$. Therefore,
$(R_i^*R_j^*-R_j^*R_i^*)(\uR ^\alpha )^*h=
(T_i^*T_j^*-T_j^*T_i^*)(\uT ^\alpha )^*h$ for $h\in \CH$,
$1\leq i, j\leq n,$ and $\alpha \in {\tilde \Lambda}$. Now
the first part of the result is clear from Proposition 4.
Further for $h\in \CL ^c(\uR)$, $R_i^*h=(R_i^c)^*h$ and so for
$h\in \CH ^c(\uT)=\CL ^c(\uR)\bigcap \CH$, $(R_i^c)^*h=R_i^*h=T_i^*h= (T_i^c)^*h.$
This proves the claim.
\qed

\begin{Definition} {\em
Let $\uT =(T_1, \ldots , T_n)$ be a contractive tuple on a Hilbert
space $\CH .$ The operator $\Delta _{\uT}=
[I-(T_1T_1^*+\cdots +T_nT_n^*)]^{\frac{1}{2}}$ is called the
{\em defect operator} of $\uT$ and
the subspace $\overline {\Delta _\uT(\CH)}$ is called the {\em defect
space\/} of $\uT .$
The tuple $\uT $ is said to be {\em pure } if
$\sum _{\alpha \in \Lambda ^m}\uT^{\alpha }(\uT^{\alpha})^*$ converges
to zero in strong operator topology as $m$ tends to infinity.}

\end{Definition}

Suppose $\sum T_iT_i^*=I$, then it is easy to see
that $\sum _{\alpha \in \Lambda ^m}\uT^{\alpha }(\uT^{\alpha})^*
=I$ for all $m$ and there is no-way this sequence can converge to zero.
So in the pure case the defect operator and the defect spaces are non-trivial.

First we restrict our attention to pure tuples. The reason for this
is that it is very easy to write down  standard dilations
  for pure tuples. So let $\CH $ be a complex, separable Hilbert space
and let $\uT $ be a pure contractive tuple on $\CH .$
Take $\tilde {\CH}=\Gamma (\mathbb{C}^n)\otimes \overline{\Delta _{\uT}(\CH)},$
and
define an operator $A:\CH \to \tilde {\CH}$ by
\begin{equation}
Ah= \sum _{\alpha }e^{\alpha }\otimes \Delta _{\uT}(\uT ^{\alpha })^*h,
\end{equation}
where the sum is taken over all $\alpha \in \tilde {\Lambda}.$
It is well-known (\cite{Po4}, \cite{AP1}) and also easily verifiable using the pureness of $\uT $, that
$A$ is an isometry with
$$A^*(e^\alpha \otimes h)=\uT ^\alpha\Delta _{\uT}h ~~ \mbox {for} ~~
\alpha \in {\tilde \Lambda}, h\in \overline {\Delta _\uT(\CH)}.$$
Now $\CH $ is considered as a subspace of
$\tilde {\CH }$ by identifying vectors $h\in \CH$ with $Ah\in {\tilde \CH}$.
Then by noting that each $V_i^*\otimes I$ leaves the range of $A$ invariant
and $\uT ^\alpha =A^*(\uV ^\alpha \otimes I)A$ for all $\alpha \in {\tilde
\Lambda}$ it is seen that the tuple ${\tilde \uV}= (V_1\otimes I, \ldots , V_n\otimes I)$
of operators on ${\tilde \CH }$ is a realization of the minimal
isometric dilation of $\uT .$ Now if $\uT$ is a commuting tuple, it is easy
to see that the range of $A$ is
 contained in ${\tilde \CH}_s=
  \Gamma _s({\mathbb{C}}^n)\otimes \overline {\Delta _{\uT}(\CH)}.$
In other words now $\CH $ can be considered as a subspace
of ${\tilde \CH}_s$. Moreover,   $\tilde {\uS}=(S_1\otimes I, \ldots , S_n\otimes I)$, as a tuple
of operators in ${\tilde \CH}_s$ is a realization of the standard commuting  dilation of
$(T_1, \ldots T_n).$ More abstractly, if $\uT$ is commuting and pure,
the standard commuting dilation of it is got by embedding $\CH$ isometrically in
$\Gamma _s(\mathbb{C}^n)\otimes \CK$, for some Hilbert space $\CK$,
such that $(S_1\otimes I_{\CK}, \ldots , S_n\otimes I_{\CK})$ is
a dilation of $\uT$ and $\overline {\mbox {~span~}}\{(\uS ^\alpha\otimes I_{\CK})h:
h\in \CH, \alpha \in \tilde {\Lambda }\}=\Gamma _s(\mathbb{C}^n)\otimes \CK$.
Up to unitary equivalence there is unique such dilation and
dim$(\CK)=$ rank$~(\Delta _{\uT}).$

\begin{Theorem}
Let $\uT $ be a pure contractive tuple on a Hilbert space $\CH$. Then the
maximal commuting piece ${\tilde \uV}^c$ of the minimal isometric dilation
${\tilde \uV}$ of $\uT$ is a realization of the
 standard commuting dilation of $\uT ^c$ if and only
if $\overline {\Delta _\uT(\CH)}=\overline {\Delta _\uT(\CH ^c(\uT))}.$
In such a case rank $(\Delta _{\uT})=$ rank $(\Delta _{\uT^c})=$
rank $(\Delta _{{\tilde \uV}})=$ rank $(\Delta _{{\tilde \uV}^c}).$
\end{Theorem}

\noindent {\sc Proof:} We denote $\CH ^c(\uT ),  \overline {\Delta _\uT(\CH)},
\overline {\Delta _\uT(\CH ^c(\uT))}$ by  $\CH^c, \CM,$ and $ \CM ^c $ respectively.
It is obvious that $\uT ^c$ is also a pure contractive
tuple. We already know from Proposition 7 that  $\tilde {\uV}^c=
(\uS\otimes I_{\CM})$ on $\Gamma _s(\mathbb{C}^n)\otimes \CM$ is
a dilation of $\uT ^c .$ It is the standard dilation if and only
if
$\CL :=\overline {\mbox {span}}\{ (\uS ^\alpha \otimes I_{\CM })Ah:
h\in \CH ^c, \alpha \in \tilde {\Lambda }\} $ is equal to
$\Gamma _s(\mathbb{C}^n)\otimes \CM $, where $A: \CH \to \tilde {\CH}$ is the isometry
defined by (2.2).

 From the definition  of $A$, using the commutativity
of the operators $T_i$,  it is clear that for $h\in \CH ^c$,
$Ah\in \Gamma _s(\mathbb{C}^n)\otimes \CM ^c.$ Hence
$\CL \subseteq \Gamma _s(\mathbb{C}^n)\otimes \CM ^c.$
Further, as $(\uS \otimes I_{\CM})$ is a dilation,
$(S_i^*\otimes I_{\CM})$ leaves $A(\CH ^c)$ invariant. Therefore,
$((I-\sum S_iS_i^*)\otimes I_{\CM})Ah\in \CL$ for $h\in \CH ^c.$ But
$(I-\sum S_iS_i^*)$ being the projection onto the vacuum space,
$( (I-\sum S_iS_i^*)\otimes I_{\CM})Ah=\omega \otimes \Delta _{\uT}h$.
As $\{\uS ^\alpha,
\alpha \in \tilde {\Lambda }\}$ spans whole of $\Gamma _s(\mathbb{C}^n)$
we get that $\Gamma _s(\mathbb{C}^n)\otimes \CM ^c \subseteq \CL$.
Hence $\CL =\Gamma _s(\mathbb{C}^n)\otimes \CM ^c$ and this way
we have proved the first claim.

Now suppose ${\tilde \uV}^c$ is a realization of the standard commuting dilation of $\uT ^c$.
This in particular means that rank$~(\Delta _{\uT ^c})=$ rank $~
(\Delta _{\tilde {\uV}^c}).$
Also as $\tilde {\uV}$ is the minimal isometric dilation of
$\uT$, rank$~(\Delta _{\uT})=$ rank$~(\Delta _{\tilde {\uV}}).$
Further as $\tilde {\uV}^c=
(\uS\otimes I_{\CM})$, rank$~(\Delta _{\tilde {\uV}^c})=$
dim$(\CM)=$ rank$~(\Delta _{\uT})$.
\qed

We may ask whether the equality of ranks in this Theorem is good enough
to make a converse statement. To answer this we make use of the following
simple lemma.

\begin{Lemma}
Suppose
$$M=\left[\begin{array}{cc}
A&B^*\\ B &C\end{array}\right]$$
is a bounded positive operator on some Hilbert space. Then rank$~(A)=$ rank
$~(\left[\begin{array}{c}
A\\ B \end{array}\right])$
\end{Lemma}

\noindent {\sc Proof:} Without loss of generality we can assume that
$M$ is a contraction. Then it is a folklore theorem that there
exists  a contraction $D$ such that $B=C^{\frac{1}{2}}DA^{\frac{1}{2}}.$
Now
$$~\left[\begin{array}{c}
A\\ B \end{array}\right]=
\left[\begin{array}{c}
A^{\frac{1}{2}}\\ C^{\frac{1}{2}}D \end{array}\right]
\left[\begin{array}{c}
A^{\frac{1}{2}} \end{array}\right],$$
and hence rank$~(\left[\begin{array}{c}
A\\ B \end{array}\right])\leq $ rank$~(A^{\frac{1}{2}}).$ But
$A$ being positive, rank $(A)$ = rank $(A^{\frac{1}{2}}).$
Therefore  rank$~(\left[\begin{array}{c}
A\\ B \end{array}\right])\leq $ rank $~(A)\leq $ rank$~(\left[\begin{array}{c}
A\\ B \end{array}\right]).$
\qed

\begin{Remark}
Let $\uT $ be a pure contractive tuple on a Hilbert space $\CH$ with
 minimal isometric dilation ${\tilde \uV}$. If rank $\Delta _{\uT}$
 and rank $\Delta _{\uT ^c}$ are finite and equal then ${\tilde \uV}^c$
 is a realization of the  standard commuting  dilation of $\uT ^c$.
\end{Remark}

\noindent {\sc Proof:} In view of Theorem 9 we need to show
that $\overline {\Delta _\uT(\CH)}=\overline {\Delta _\uT(\CH ^c(\uT))}.$
Since $\overline {\Delta _\uT(\CH)}\supseteq\overline {\Delta _\uT(\CH ^c(\uT))},$
and these spaces are now finite dimensional, it suffices to show that
their dimensions are equal or  rank$~(\Delta _{\uT})=$ rank$~(\Delta _{\uT}P_{\CH ^c}).$
Clearly rank$~(\Delta _{\uT})\geq $ rank$~(\Delta _{\uT}P_{\CH ^c}).$
Also by assumption, rank$~(\Delta _{\uT})=$ rank$~(\Delta _{\uT ^c}).$
By positivity rank$~(\Delta _{\uT ^c})=$
rank$~(\Delta _{\uT^c}^2).$ And then by previous Lemma
rank$~(\Delta _{\uT^c}^2)=$ rank$~(P_{\CH ^c}(\Delta _{\uT}^2)P_{\CH ^c})=$
rank$~(\Delta _{\uT}^2P_{\CH ^c})\leq $
rank$~(\Delta _{\uT}P_{\CH ^c}).$
\qed

If both the ranks are infinite then we can not ensure that
$\overline {\Delta _\uT(\CH)}=\overline {\Delta _\uT(\CH ^c(\uT))}$
is seen by the following example.
\begin{Example}{\em
Let $\uR=(R_1, R_2)$ be a commuting pure contractive 2-tuple on
an infinite dimensional Hilbert space $\CH _0$ (We can even take $R_1, R_2$
as scalars) such that $\overline {\Delta _\uR(\CH _0)}$ is
infinite dimensional.
 Take $\CH =\CH _0\oplus \mathbb{C}^2,$ and let $T_1, T_2$ be
operators on $\CH $ defined by
$$T_1=\left[\begin{array}{ccc}
R_1 && \\
&0& t_1\\
&0 & 0\end{array}\right],
T_2=\left[\begin{array}{ccc}
R_2 && \\
&0& 0\\
&t_2& 0\end{array}\right],$$
where $t_1, t_2$ are any two scalars, $0<t_1, t_2<1.$ Then $\uT=(T_1, T_2)$
is a pure contractive tuple. Making use of Corollary 5,
 $\CH ^c(\uT)=
\CH _0$ (thought of as a subspace of $\CH $ in the natural way) and
the maximal commuting piece of $\uT$ is $(R_1, R_2)$, and therefore
rank$~(\Delta _{\uT^c})=$ rank$~(\Delta _{\uT})=\infty .$ But
$\overline {\Delta _\uT(\CH)}=\overline {\Delta _\uR(\CH _0)}\oplus
\mathbb{C}^2$.}
\end{Example}

We do not know how to extend Theorem 9 to contractive tuples which
are not necessarily pure.

\end{section}

\begin{section}{Commuting Tuples}
\setcounter{equation}{0}

In this Section we wish to consider commutative contractive tuples.
 Let us begin with describing the way one obtains
two standard dilations for such tuples.

 Recall  standard tuples $\uV$ and $\uS$ on
Fock spaces $\Gamma ({\mathbb{C}}^n)$, and $
\Gamma _s({\mathbb{C}}^n)$ respectively, introduced in the Introduction.
Let
$C^*(\uV)$, and  $C^*(\uS)$ be unital $C^*$ algebras generated by them.
For any $\alpha , \beta \in {\tilde \Lambda }$, $\uV ^{\alpha }(I-\sum V_iV_i^*)
(\uV ^{\beta })^*$ is the rank one operator $x\mapsto  \langle  e^{\beta }, x\rangle e^\alpha ,$
formed by basis vectors $e^{\alpha }, e^{\beta }$.  So $C^*(\uV)$ contains
all compact operators. In a similar way we see that $C^*(\uS)$ also contains
all compact operators of $\Gamma _s({\mathbb{C}}^n)$.
As $V_i^*V_j=\delta _{ij},$ it is easy to see that $C^*(\uV)=
~~\overline {\mbox {span}}~~\{\uV^\alpha (\uV^\beta)^*: \alpha , \beta \in
{\tilde \Lambda}\}.$  By explicit computation  commutators $[S_i^*, S_j]$
are compact for all $i,j$ (See  \cite{Ar2}, Proposition 5.3, or \cite{BB}). Therefore we can also obtain
$C^*(\uS)=~~\overline {\mbox {span}}~~\{\uS^\alpha (\uS^\beta)^*: \alpha , \beta \in
{\tilde \Lambda}\}.$

Suppose $\uT $ is a contractive tuple on a Hilbert space $\CH$. We obtain
a certain completely positive map (Popescu's Poisson transform) from $C^*(\uV)$ to $\CB (\CH)$,  as follows.
For $0< r <1 $ the tuple $r\uT= (rT_1, \ldots , rT_n)$ is clearly a pure contraction.
So by (2.2) we have an isometry $A_r: \CH \to \Gamma ({\mathbb{C}}^n)\otimes
\overline {\Delta _r(\CH)}$ defined by
$$A_rh = \sum _{\alpha }e^{\alpha }\otimes \Delta _r((r\uT) ^\alpha )^*h, ~~h\in \CH ,$$
where $\Delta _r=(I-r^2\sum T_iT_i^*)^{\frac{1}{2}}.$ So for every $0<r<1$
we have a completely positive map $\psi _r: C^*(\uV)\to \CB(\CH)$
defined by
$$\psi _r(X)=A_r^*(X\otimes I)A_r, ~~X\in C^*(\uV).$$
By taking limit as $r$ increases to 1 (See \cite{Po4} or \cite{AP1} for details), we obtain a
unital completely positive map $\psi $ from $C^*(\uV)$ to $\CB(\CH)$
satisfying
$$\psi (\uV^\alpha (\uV^\beta)^*)= \uT^\alpha (\uT^\beta )^* ~~\mbox{for} ~ \alpha , \beta \in
{\tilde \Lambda}.$$
As $C^*(\uV)=
~~\overline {\mbox {span}}~~\{V^\alpha (V^\beta)^*: \alpha , \beta \in
{\tilde \Lambda}\},$ $\psi $ is the unique such completely positive map.
Now consider the minimal Stinespring dilation of $\psi .$
So we have a Hilbert space ${\tilde \CH}$ containing $\CH $,  and a
unital $*$-homomorphism $\pi : C^*(\uV)\to \CB({\tilde \CH }),$ such that
$$\psi (X)=P_{\CH}\pi(X)|_{\CH} ~~\forall X\in C^*(\uV),$$
and $\overline {\mbox {span}}~\{ \pi(X)h: X\in C^*(\uV), h\in \CH\}=
{\tilde \CH}.$
Taking ${\tilde \uV}=({\tilde V_1}, \ldots , {\tilde V_n})
= (\pi(V_1), \ldots , \pi(V_n))$, one verifies that each $\tilde {(V_i)}^*$
leaves $\CH $ invariant and ${\tilde \uV}$ is the unique minimal isometric
dilation of $\uV .$

In a similar fashion if $\uT$ is commuting by considering $C^*(\uS)$ instead
of $C^*(\uV)$, and restricting $A_r$ in the range to $\Gamma _s({\mathbb{C}}^n)$,
and taking limits as before (See \cite{Ar2}, \cite{Po4}, \cite{AP1}) we obtain the unique unital
completely positive map $\phi : C^*(\uS)\to \CB (\CH)$, satisfying
$$\phi (\uS^\alpha (\uS^\beta)^*)= \uT^\alpha (\uT^\beta )^* ~~~~~ \alpha , \beta \in
{\tilde \Lambda}.$$
Consider the minimal Stinespring dilation of $\phi .$ Here we have a Hilbert
space $\CH _1 $ containing  $\CH$ and a unital $*$-homomorphism
$\pi _1 : C^*(\uS)\to \CB(\CH _1),$ such that
$$\phi (X)=P_{\CH }\pi _1(X)|_{\CH} ~~~~\forall X\in C^*(\uS),$$
and $\overline {\mbox {span}}~\{ \pi _1(X)h: X\in C^*(\uS), h\in \CH \}=
{ \CH _1}.$
Taking ${\tilde \uS}=({\tilde S_1}, \ldots , {\tilde S_n})
= (\pi _1(S_1), \ldots , \pi _1(S_n))$,  $\tilde {\uS}$ is the standard
commuting dilation of $\uT $ by definition (It is not difficult to verify
that it is a minimal dilation in the sense of our Definition 2). As minimal
Stinespring dilation is unique up to unitary equivalence, standard commuting
dilation is also unique up to unitary equivalence.

\begin{Theorem} (Main Theorem)
Suppose $\uT$ is a commuting contractive tuple on a Hilbert space $\CH .$ Then
  the maximal commuting piece of
the minimal isometric dilation of $\uT $ is a realization of the
standard commuting dilation of $\uT$.
\end{Theorem}

Our approach to prove this theorem is as follows. First we consider the standard
commuting dilation of $\uT $ on a Hilbert space $\CH _1$ as described above.
Now the standard tuple $\uS $ is also a contractive tuple. So we have a unique
unital completely positive map $\eta : C^*(\uV)\to C^*(\uS)$, satisfying
$$\eta (\uV^\alpha (\uV^\beta)^*)= \uS^\alpha (\uS^\beta )^* ~~~~ \alpha , \beta \in
{\tilde \Lambda}.$$
Now clearly $\psi = \phi \circ \eta $. Consider the minimal Stinespring dilation
of the composed map $\pi _1\circ \eta : C^*(\uV)\to \CB(\CH _1)$. Here we obtain a Hilbert space
 $\CH _2 $ containing  $\CH _1$ and a unital $*$-homomorphism
$\pi _2 : C^*(\uV)\to \CB(\CH _2),$ such that
$$\pi _1\circ \eta (X)=P_{\CH _1 }\pi _2(X)|_{\CH _1}, ~~ ~~\forall X\in C^*(\uV),$$
and $\overline {\mbox {span}}~\{ \pi _2(X)h: X\in C^*(\uV), h\in \CH _1\}=
{\CH _2}.$
Now we have a commuting diagram as follows

\hskip1.5in \begin{picture}(200,115)

\put(0,20){$C^*(\uV)$}

\put(40,21){$\longrightarrow$}

\put(70,20){$C^*(\uS)$}

\put(110,21){$\longrightarrow$}

\put(140,20){$\mathcal{B}(\CH)$}

\put(140,60){$\mathcal{B}(\CH_1)$}

\put(140,100){$\mathcal{B} (\CH_2)$}

\put(40,30){\vector(4,3){80}}

\put(110,30){\vector(4,3){20}}

\put(150,40){$\downarrow $}

\put(150,80){$\downarrow $}

\put(45,10){$\eta$}

\put(115,10){$\phi$}

\put(110,40){$\pi_1$}

\put(70,70){$\pi_2$}

\end{picture}

\noindent where all the down arrows are
compression maps, horizontal arrows are unital completely positive maps and
diagonal arrows are unital $*$-homomorphisms.

Taking ${\hat  \uV}=({\hat V_1}, \ldots , {\hat V_n})
= (\pi _2(V_1), \ldots , \pi _2(V_n))$, we need to show (i) ${\hat \uV}$ is
the minimal isometric dilation of $\uT$ and (ii) ${\tilde {\uS }}
=(\pi _1(S_1), \ldots , \pi _1(S_n))$ is the
maximal commuting piece of $\hat {\uV}$. Due to uniqueness up to unitary equivalence
of minimal Stinespring dilation, we have (i) if we can show that $\pi _2$ is
a minimal dilation of $\psi=\phi \circ\eta .$ For proving this we actually make
use of (ii). At first  we prove (ii) in a very special case.
\begin{Definition}{\em
A $n$-tuple $\uT =(T_1, \ldots , T_n) $ of operators on a Hilbert space $\CH $
is called a spherical unitary if it is commuting, each $T_i$ is normal,
and $T_1T_1^*+\cdots +T_nT_n^*=I.$}
\end{Definition}

Actually, if $\CH $ is a finite dimensional Hilbert space and $\uT $
is a commuting tuple on $\CH $ satisfying $\sum T_iT_i^* =I$, then
it is automatically a spherical unitary, that is, each $T_i$ is normal.
This is the case because here standard commuting
dilation of $\uT $ is a tuple of normal operators  and
hence each $T_i^*$ is subnormal (or see \cite{At1} for this result) and
all finite dimensional subnormal
operators are normal (see \cite{Ha}).

Note that  if $\uT $ is a spherical unitary we have $\phi (\uS^\alpha
(I-\sum S_iS_i^*)(\uS ^\beta )^*)
=\uT ^\alpha (I-\sum T_iT_i^*)(\uT ^\beta)^*=0$ for any
$\alpha , \beta \in {\tilde \Lambda}.$ This forces that
$\phi (X)=0$ for any compact operator $X$ in $C^*(\uS).$
Now as the commutators $[S_i^*, S_j]$ are all compact we
see that $\phi $ is a unital $*$-homomorphism. So the minimal
Stinespring dilation of $\phi $ is itself. So the following result
yields
Theorem 13 for spherical unitaries.
\begin{Theorem}
Let $\uT$ be a spherical unitary on a Hilbert space $\CH .$ Then
the maximal commuting piece of the minimal isometric dilation of
$\uT $ is $\uT $.
\end{Theorem}

As proof of this Theorem involves some lengthy computations
we prefer to
postpone it. But assuming this, we prove the main Theorem.

\noindent {\sc Proof of Theorem 13 :} As $C^*(\uS)$ contains the ideal of all compact
operators by standard $C^*$-algebra theory we have a direct
sum decomposition of $\pi _1$ as follows. Take
${\mathcal H}_1 = {\mathcal H}_{1C} \oplus {\mathcal H}_{1N}$\\
where ${\mathcal H}_{1C} = \overline{\mbox{span}}\{\pi(X)h : h\in {\mathcal H},
X \in
C^*(\underline{S})$ and $X$ is compact$\}$ and
 ${\mathcal H}_{1N} = {\mathcal H}_1 \ominus {\mathcal H}_{1C}$,
Clearly ${\mathcal H}_{1C}$ is a reducing subspace for $\pi _1$.
Therefore \[\pi _1(X)=\begin{pmatrix}
  \pi_{1C}(X)&  \\
     & \pi_{1N}(X)
\end{pmatrix}     \]
that is, $\pi _1=\pi _{1C}\oplus \pi _{1N}$ where $\pi_{1C}(X) = P_{{\mathcal H}_{1C}} \pi _1(X) P_{{\mathcal H}_{1C}}$,
$\pi_{1N}(X) = P_{{\mathcal H}_{1N}} \pi _1(X) P_{{\mathcal H}_{1N}}$.
As observed by Arveson \cite{Ar2}, $\pi _{1C}(X)$ is just the identity representation
with some multiplicity. More precisely, $\CH _{1C}$ can be factored
as $\CH _{1C}= \Gamma _s({\mathbb{C}}^n)\otimes \overline {\Delta _{\uT}(\CH)},$
such that $\pi _{1C}(X)=X\otimes I$, in particular
$\pi _{1C}(S_i)=S_i\otimes I$. Also $\pi _{1N}(X)=0$ for compact
$X$. Therefore, taking $Z_i=\pi _{1N}(S_i)$, $\underline {Z}=(Z_1, \ldots , Z_n)$
is a spherical unitary.

Now as $\pi _1\circ \eta =(\pi _{1C}\circ \eta )\oplus (\pi _{1N}\circ \eta)$
and
the minimal Stinespring dilation of a direct sum of two completely
positive maps is the direct sum of minimal Stinespring dilations. So $\CH _2
$  decomposes as $\CH _2= \CH _{2C}\oplus \CH _{2N},$ where
$\CH _{2C}, \CH _{2N}$ are orthogonal reducing subspaces of $\pi _2$,
such that $\pi _2$ also decomposes, say $\pi _2=\pi _{2C}\oplus
\pi _{2N},$ with
$$\pi _{1C}\circ \eta (X)= P_{\CH _{1C}}\pi _{2C}(X)|_{\CH _{1C}}, ~~
\pi _{1N}\circ \eta (X)= P_{\CH _{1N}}\pi _{2N}(X)|_{\CH _{1N}},$$
for $
X\in C^*(\uV)$
 with
$\CH _{2C}= ~\overline {\mbox {span}}~\{
\pi _{2C}(X)h: X\in C^*(\uV), h\in \CH _{1C}\}$ and
$\CH _{2N}= ~\overline {\mbox {span}}~\{
\pi _{2N}(X)h: X\in C^*(\uV), h\in \CH _{1N}\}$.  It is also
not difficult to see that $\CH _{2C}= ~\overline {\mbox {span}}~\{
\pi _{2C}(X)h: X\in C^*(\uV), X ~~\mbox {compact},~  h\in \CH _{1C}\}$ and
hence $\CH _{2C}$ factors as $\CH _{2C}= \Gamma ({\mathbb{C}}^n)\otimes
\overline {\Delta _{\uT}(\CH)}$ with $\pi _{2C}(V_i)=V_i\otimes I.$
Also $(\pi _{2N}(V_1), \ldots , \pi _{2N}(V_n))$ is a minimal
isometric dilation of spherical isometry $(Z_1, \ldots , Z_n)$.
Now by  Proposition 6, Theorem 15 and Corollary 5, we get that
$(\pi _1(S_1), \ldots , \pi _1(S_n))$ acting on $\CH _1$ is the
maximal commuting piece of $(\pi _2(V_1), \ldots , \pi _2(V_n))$.

All that remains to show is that $\pi _2$ is the minimal Stinespring
dilation of $\phi \circ \eta .$ Suppose this is not the case.
Then we get a reducing subspace $\CH _{20}$ for $\pi _2$ by
taking
$\CH _{20} =~\overline {\mbox {span}}~\{\pi _2(X)h: X\in C^*(\uV), h\in \CH
\} .$
Take $\CH _{21}=\CH _2\ominus \CH _{20}$ and correspondingly
decompose $\pi _2$ as $\pi _2= \pi _{20}\oplus \pi _{21}$,
\[\pi _2(X)=\begin{pmatrix}
  \pi_{20}(X)&  \\
     & \pi_{21}(X)
\end{pmatrix}     \]

 Note
that we already have $\CH \subseteq \CH _{20}.$  We claim that
$\CH _2\subseteq \CH _{20}.$
Firstly, as $\CH _1$ is the space where the maximal commuting piece
of $(\pi _2(V_1), \ldots , \pi _2(V_n))= (\pi _{20}(V_1)\oplus
\pi _{21}(V_1), \ldots , \pi _{20}(V_n)\oplus \pi _{21}(V_n))$ acts,
by the first part of Corollary 5, $\CH _1$ decomposes as
$\CH _1=\CH _{10}\oplus \CH _{11}$ for some subspaces $\CH _{10}
\subseteq \CH _{20},$ and
$\CH _{11}\subseteq  \CH _{21}$.  So for
$X\in C^*(\uV)$, $P_{\CH _1}\pi _2(X)P_{\CH _1}$, has the form (see the
diagram)
$$ P_{\CH _1}\pi _2(X)P_{\CH _1}=
\begin{pmatrix}
  \pi _{10}\circ \eta (X)&0& & \\
  0 & 0& & \\
   & & \pi_{11}\circ \eta(X)& 0 \\
   & & 0 & 0
\end{pmatrix}
$$
where $\pi _{10}, \pi _{11}$ are compressions of $\pi _1$ to
$\CH _{10}$, $\CH _{11}$ respectively. As the mapping $\eta $ from $C^*(\uV)$
to $C^*(\uS)$ is clearly surjective, it follows that $\CH _{10},
\CH _{11}$ are reducing subspaces for $\pi _{1}.$ Now as $\CH $
is contained in $\CH _{20}$, in view of minimality of $\pi _1$
as a Stinespring dilation, $\CH _1\subseteq \CH _{20}$. But then
the minimality of $\pi _2$ shows that $\CH _2\subseteq \CH _{20}$.
Therefore, $\CH _2=\CH _{20}.$
\qed

\noindent {\sc Proof of Theorem 15 :} Here we need a different presentation of the minimal
isometric dilation. This is known as  Sch\"{a}ffer construction \cite{Sc} in the
one variable case, and \cite{Po1} is a good reference for the multivariate
case.
Here we decompose the dilation space $\tilde {\CH} $ as $ \tilde
{\CH } ={\mathcal H} \oplus
(\Gamma({\mathbb{C}}^n)\otimes\mathcal{D})$
where $\mathcal{D}$ is the closure of the range of  operator
\[ D:\underbrace{{\mathcal H}\oplus \cdots \oplus {\mathcal H}}_{n~
copies} \rightarrow \underbrace{{\mathcal H} \oplus \cdots \oplus {\mathcal H}}_{n
~copies} \]
and $D$ is the positive square root of
\[ D^2=[\delta_{ij}I-T_i^*T_j]_{n \times n}. \]
Whenever it is convenient for us we identify $\underbrace{{\mathcal H}\oplus \cdots \oplus {\mathcal
H}}_{n~copies}$ with ${\mathbb{C}}^n\otimes {\mathcal H}$ so that
\[ (h_1, \ldots ,h_n)= \sum_{i=1}^n e_i \otimes h_i. \]
Then
\begin{equation}
 D(h_1, \ldots ,h_n)= D(\sum_{i=1}^n e_i \otimes h_i)=\sum_{i=1}^n
e_i\otimes (h_i - \sum _{j=1}^nT^*_iT_j h_j).
\end{equation}
And the minimal isometric dilation $\tilde{V_i}$ has the form
\begin{equation}
 \tilde{V_i}(h\oplus \sum_{\alpha \in
{\tilde{\Lambda}}}e^{\alpha}\otimes
d_{\alpha})= T_ih\oplus D(e_i\otimes h )\oplus e_i\otimes(\sum_{\alpha \in
{\tilde{\Lambda}}}
e^{\alpha}\otimes d_{\alpha})
\end{equation}
for $h \in \CH$,
 $d_{\alpha} \in
\mathcal{D}$ for $\alpha \in \tilde{\Lambda}$, and $1 \leq i \leq n$
($\mathbb{C}^n\omega \otimes \CD$ has been identified with $\CD$).

In the present case as $\sum T_iT_i^*=I$, by direct computation
 $D^2$ is seen to be a projection. So, $D$ which is the positive
square root of $D^2$, is equal to $D^2$. Also by Fuglede-Putnam
theorem (\cite{Ha},  \cite{Pu}), $\{T_1, \ldots ,T_n,T_1^*, \ldots ,T_n^*\}$ forms a commuting family
of operators. Then we get
\begin{eqnarray}
 D(h_1, \ldots ,h_n)
&= &  \sum_{i,j=1}^n e_i \otimes T_j (T_j^*h_i - T_i^*h_j)
 = \sum_{i,j=1}^n e_i \otimes T_j (h_{ij})
\end{eqnarray}
where $h_{ij} = T_j^*h_i - T_i^*h_j$ for $1 \leq i, j \leq n$.
Note that $h_{ii}=0$ and $h_{ji}=-h_{ij}.$

Now we apply Proposition 4 to the tuple $\tilde {\uV}$ acting
on $\tilde {\CH }.$ Suppose $y\in
 \CH ^{\perp}\bigcap \tilde {\CH }^c(\tilde {\uV}). $ We wish
 to show that $y=0$. We assume $y\neq 0$ and arrive at a contradiction.
 One can decompose $y$ as $y=0\oplus \sum _{\alpha \in \tilde {\Lambda }}
 e^{\alpha }\otimes y_{\alpha },$ with $y_{\alpha }\in \CD .$ If for
 some $\alpha $, $y_{\alpha }\neq 0$, then
 $\langle \omega \otimes y_{\alpha }, (\tilde {\uV}^{\alpha })^*y\rangle
 =\langle e^{\alpha }\otimes y_{\alpha }, y\rangle =
 \langle y_{\alpha }, y_{\alpha }\rangle \neq  0.$
 Since each $(\tilde {V}_i)^*$ leaves $\tilde {\CH }^c(\tilde {\uV})$
 invariant, $(\tilde {\uV}^{\alpha })^*y\in \tilde {\CH }^c(\tilde {\uV}).$
 So without loss of generality we can assume $\|y_0\|=1.$

Taking $\tilde {y}_m =\sum _{\alpha \in \Lambda ^m}e^{\alpha }\otimes y_{\alpha}$,
we get $y=0\oplus \oplus _{m\geq 0}(\tilde {y}_m).$
 As $y_0 \in {\mathcal D}$,
 $y_0=
D(h_1, \ldots ,h_n)$, for some $(h_1, \ldots ,h_n)$ (Presently $D$ being a projection
its range is closed). Set
$\tilde{x_0}=\tilde{y_0}=y_0$, and for $m\geq 1$,
$$ \tilde{x}_m = \sum_{i_1, \ldots ,i_{m-1},i,j=1}^n e_{i_1}
\otimes \cdots \otimes e_{i_{m-1}} \otimes e_i \otimes D(e_j\otimes
T_{i_1}^*
\ldots T_{i_{m-1}}^* h_{ij}).$$
Clearly  $\tilde{x}_m \in ({\mathbb
C}^n)^{\otimes
m} \otimes {\mathcal D}$ for all $m \in
{\mathbb{N}}$.
From the definition  (3.2) of $\tilde {V}_i$, commutativity of
the operators $T_i$,    and the fact that $D$ is projection,
we have
\begin{eqnarray*}
\sum_{1\leq i < j \leq n} (\tilde{V}_i \tilde{V}_j -
\tilde{V}_j
\tilde{V}_i)h_{ij}
&=&\sum _{1\leq i<j\leq n}(T_iT_jh_{ij}-T_jT_ih_{ij})\\
 & &+\sum_{1\leq i < j \leq n}D(e_i\otimes T_j h_{ij} -
e_j\otimes T_i h_{ij})\\
& &
  +\sum_{1\leq i < j \leq n}(e_i\otimes D(e_j\otimes h_{ij})
   -e_j\otimes D(e_i \otimes h_{ij}))\\
    & =& D\{\sum_{1\leq i < j \leq n} (e_i\otimes T_j h_{ij} -
e_j\otimes T_i h_{ij})\}
  +\sum_{ i,j=1}^n e_i\otimes D(e_j\otimes h_{ij})\\
& =& D(\sum _{1\leq i,j\leq n}e_i\otimes T_jh_{ij}) + \sum_{ i,j=1}^n e_i\otimes D(e_j\otimes h_{ij})\\
& =& D^2(h_1, \ldots ,h_n) + \sum_{ i,j=1}^n e_i\otimes
D(e_j\otimes h_{ij})\\
&=& \tilde{x}_0 +\tilde{x}_1.
\end{eqnarray*}
Therefore $\langle y, \tilde{x}_0 +\tilde{x}_1 \rangle =0$ by
Proposition 4. Now for $m \geq 2$.
\begin{eqnarray*}
& & \sum_{i_1, \ldots ,i_{m-1}=1}^n \tilde{V}_{i_1} \ldots
\tilde{V}_{i_{m-1}}(\sum_{i,j=1}^n
(\tilde{V}_i\tilde{V}_j-\tilde{V}_j\tilde{V}_i)T_{i_1}^*\ldots T_{i_{m-2}}^*
T_j^*h_{i_{m-1}i})\\
& =& \sum_{i_1, \ldots ,i_{m-1}=1}^n \tilde{V}_{i_1} \ldots
\tilde{V}_{i_{m-1}}[\sum_{i,j=1}^n D(e_i\otimes
T_jT_{i_1}^*\ldots T_{i_{m-2}}^*T_j^*h_{i_{m-1}i} - e_j\otimes
T_iT_{i_1}^*\ldots T_{i_{m-2}}^*T_j^*h_{i_{m-1}i})\\
& &+\sum_{i,j=1}^n \{e_i \otimes D(e_j \otimes T_{i_1}
^*\ldots T_{i_{m-2}}^*T_j^*h_{i_{m-1}i}) - e_j\otimes
D(e_i \otimes T_{i_1}^*\ldots T_{i_{m-2}}^*T_j^*h_{i_{m-1}i})\}]
\end{eqnarray*}
\begin{eqnarray*}
&= &\sum_{i_1,\ldots ,i_{m-1}=1}^n e_{i_1} \otimes \cdots \otimes
e_{i_{m-1}} \otimes \\
& & [D(\sum_{i,j=1}^n e_i \otimes T_jT_{i_1}^* \ldots
T_{i_{m-2}}^*T_j^* h_{i_{m-1}i}
-e_j \otimes T_iT_{i_1}^* \ldots  T_{i_{m-2}}^*T_j^* h_{i_{m-1}i})\\
& & \{\sum_{i,j=1}^n e_i
\otimes D(e_j\otimes T_{i_1}^* \ldots   T_{i_{m-2}}^*T_j^* h_{i_{m-1}i})
-\sum_{i,j=1}^n e_i \otimes D(e_j\otimes T_{i_1}^*
\ldots   T_{i_{m-2}}^*T_i^* h_{i_{m-1}j})\}]\\
& & \mbox{(in the term above, $i$ and $j$ have been interchanged in the
last summation)}\\
& =& \sum_{i_1,\ldots  ,i_{m-1}=1}^n e_{i_1} \otimes \cdots
\otimes e_{i_{m-1}} \otimes \\
& & [D(\sum _{i=1}^n e_i \otimes T_{i_1}^* \ldots  T_{i_{m-2}}^* h_{i_{m-1}i}
- \sum _{i, j=1}^ne_j
\otimes T_iT_{i_1}^* \ldots   T_{i_{m-2}}^*T_j^* h_{i_{m-1}i} )\\
& & \{\sum_{i,j=1}^n e_i \otimes D(e_j\otimes
(T_{i_1}^* \ldots   T_{i_{m-2}}^*T_j^* h_{i_{m-1}i}
- T_{i_1}^* \ldots   T_{i_{m-2}}^*T_i^* h_{i_{m-1}j})\}]\\
& =& \sum_{i_1,\ldots  ,i_{m-1}=1}^n e_{i_1} \otimes \cdots \otimes
e_{i_{m-1}} \otimes \sum_{i=1}^n
D(e_i\otimes T_{i_1}^* \ldots   T_{i_{m-2}}^* h_{i_{m-1}i})\\
& & +\sum_{i_1,\ldots  ,i_{m-1},i,j=1}^n e_{i_1} \otimes \cdots  \otimes
e_{i_{m-1}} \otimes e_i \otimes \\
& & D(e_j\otimes T_{i_1}^* \ldots
T_{i_{m-2}}^*(T_j^*T_{i}^*h_{i_{m-1}} - T_j^*T_{i_{m-1}}^*h_{i} -
 T_i^*T_{j}^*h_{i_{m-1}}+T_i^*T_{i_{m-1}}^*h_j))\\
&= &   \sum_{i_1,\ldots  ,i_{m-2},i,j=1}^n e_{i_1} \otimes \cdots  \otimes
e_{i_{m-2}} \otimes ( e_i
\otimes D(e_j \otimes T_{i_1}^* \ldots   T_{i_{m-2}}^* h_{ij})\\
&  & +\sum_{i_1,\ldots  ,i_{m-1},i,j=1}^n e_{i_1} \otimes \cdots  \otimes
e_{i_{m-1}} \otimes e_i \otimes D(e_j\otimes T_{i_1}^* \ldots
T_{i_{m-2}}^*( - T_j^*T_{i_{m-1}}^*h_{i}+T_i^*T_{i_{m-1}}^*h_j))\\
& & \mbox{(in the  term above, index $i_{m-1}$ has been replaced
by $i$ and $i$ has been replaced by} \\
& & \mbox{ $j$ in the first summation)}\\
& =& \sum_{i_1,\ldots  ,i_{m-2},i,j=1}^n e_{i_1} \otimes \cdots \otimes
e_{i_{m-2}} \otimes e_i
\otimes D(e_j \otimes T_{i_1}^* \ldots   T_{i_{m-2}}^* h_{ij})\\
& & -\sum_{i_1,\ldots  ,i_{m-1},i,j=1}^n e_{i_1} \otimes \cdots \otimes
e_{i_{m-1}} \otimes e_i
\otimes D(e_j\otimes T_{i_1}^* \ldots   T_{i_{m-1}}^* h_{ij})
=\tilde{x}_{m-1} - \tilde{x}_m.
\end{eqnarray*}
So, $\langle y, \tilde{x}_{m-1} - \tilde{x}_m\rangle =0.$

Next, we would show that $\|\tilde{x}_{m+1}\| = \| \tilde{x}_0 \| =1$
for all $m \in {\mathbb{N}}$.
\begin{eqnarray*}
\|\tilde{x}_{m+1}\|^2 & =& \langle
\sum_{i_1,\ldots  ,i_m,i,j=1}^n
e_{i_1} \otimes \cdots  \otimes e_{i_m} \otimes e_i
\otimes D(e_j\otimes T_{i_1}^* \ldots   T_{i_m}^* h_{ij}),\\
& & \sum_{i_1^{\prime},\ldots  ,i_{m-1}^{\prime},i^{\prime},j^{\prime}=1}^n
e_{i_1^{\prime}} \otimes \cdots \otimes
e_{i_m^{\prime}} \otimes e_{i^{\prime}}
\otimes D(e_{j^{\prime}}\otimes T_{i_1^{\prime}}^*
\ldots   T_{i_m^{\prime}}^* h_{i^{\prime}j^{\prime}})\rangle\\
&= &  \sum_{i_1,\ldots  ,i_m,i=1}^n \langle \sum_{j=1}^n D(e_j\otimes
T_{i_1}^* \ldots   T_{i_m}^* h_{ij}), \sum_{j^{\prime}=1}^n
D(e_{j^{\prime}}\otimes
T_{i_1}^* \ldots   T_{i_m}^*
h_{ij^{\prime}})\rangle\\
&= &  \sum_{i_1,\ldots  ,i_m,i=1}^n \langle D(\sum_{j=1}^n e_j\otimes
T_{i_1}^* \ldots   T_{i_m}^* h_{ij}),
\sum_{j^{\prime}=1}^n e_{j^{\prime}}\otimes T_{i_1}^* \ldots   T_{i_m}^*
 h_{ij^{\prime}}\rangle\\
& =& \sum_{i_1,..,i_m,i=1}^n \langle \sum_{l,k=1}^n e_l\otimes
T_k (T_k^*T_{i_1}^* \ldots   T_{i_m}^* h_{il}- T_l^* T_{i_1}^*\ldots
T_{i_m}^* h_{ik}),\sum_{j^{\prime}=1}^n e_{j^{\prime}}\otimes T_{i_1}^*
\ldots   T_{i_m}^* h_{ij^{\prime}}\rangle\\
& =& \sum_{i_1,..,i_m,i,j=1}^n \langle \sum_{k=1}^n T_k
(T_k^*T_{i_1}^* \ldots   T_{i_m}^* h_{ij}- T_j^* T_{i_1}^*\ldots   T_{i_m}^*
h_{ik}), T_{i_1}^* \ldots   T_{i_m}^* h_{ij}\rangle\\
& =& \sum_{i,j=1}^n\langle \sum_{k=1}^n (T_kT_k^*h_{ij}-
T_kT_j^*h_{ik}),h_{ij}\rangle\\
& =& \sum_{i,j=1}^n\langle \sum_{k=1}^n(T_kT_k^*T_j^*h_i -
T_kT_k^*T_i^*h_j - T_kT_j^*T_k^*h_i + T_kT_j^*T_i^*h_k),
T_j^*h_i - T_i^*h_j \rangle\\
& =& \sum_{i,j=1}^n\langle \sum_{k=1}^n(T_kT_j^*T_i^*h_k -
T_kT_k^*T_i^*h_j), T_j^*h_i - T_i^*h_j \rangle\\
& =& \sum_{i,j=1}^n\langle \sum_{k=1}^n(T_kT_j^*T_i^*h_k) -
T_i^*h_j, T_j^*h_i - T_i^*h_j \rangle\\
& =&\sum_{i,j=1}^n \langle T_j (\sum_{k=1}^nT_kT_j^*T_i^*h_k) -
T_jT_i^*h_j, h_i \rangle
-\sum_{i,j=1}^n\langle T_i (\sum_{k=1}^nT_kT_j^*T_i^*h_k) -
T_iT_i^*h_j, h_j \rangle\\
& =& \sum_{i=1}^n \langle
\sum_{j=1}^nT_jT_j^*(\sum_{k=1}^nT_kT_i^*h_k)
- \sum_{j=1}^n T_jT_i^*h_j, h_i \rangle
-\\
& & \sum_{j=1}^n \langle
\sum_{i=1}^nT_iT_i^*(\sum_{k=1}^nT_kT_j^*h_k)
- \sum_{i=1}^n T_iT_i^*h_j, h_j \rangle\\
& =& \sum_{i=1}^n \langle (\sum_{k=1}^nT_kT_i^*h_k) -
\sum_{j=1}^n T_jT_i^*h_j, h_i \rangle
-\sum_{j=1}^n \langle \sum_{k=1}^nT_kT_j^*h_k- h_j, h_j \rangle\\
& =& \sum_{j=1}^n \langle h_j -\sum_{k=1}^n T_kT_j^*h_k), h_j \rangle
=\langle D(h_1,\ldots  ,h_n),(h_1,\ldots  ,h_n) \rangle
=\| \tilde{x}_0 \| ^2 =1.
\end{eqnarray*}
As $\langle y, \tilde{x}_0 +\tilde{x}_1\rangle =0$ and $\langle y,
\tilde{x}_m -\tilde{x}_{m+1} \rangle =0
$ for $m \in {\mathbb{N}}$,
 we get $\langle y, \tilde{x}_0 +\tilde{x}_{m+1} \rangle =0$
for $m \in {\mathbb{N}}$.  This implies
$1=\langle \tilde {y}_0, \tilde {y_0}\rangle =
\langle \tilde{y}_0,\tilde{x}_0 \rangle = - \langle
\tilde{y}_{m+1}, \tilde{x}_{m+1} \rangle$. By Cauchy-Schwarz
inequality,
$1 \leq \| \tilde{y}_{m+1} \| \| \tilde{x}_{m+1} \|$ , i.\@ e., $1
\leq \| \tilde{y}_{m+1} \|$ for $m \in {\mathbb{N}}$. This is a
contradiction as $y=0\oplus \oplus _{m\geq 0}\tilde {y}_m$ is in
the  Hilbert space $\tilde {\CH }$.
\qed

\end{section}

\begin{section}{Representations of Cuntz Algebras}
\setcounter{equation}{0}

For $n\geq 2$, the Cuntz algebra $\CO _n$ is the $C^*$-algebra
generated by $n$-isometries $\us= \{s_1, \ldots , s_n\}$, satisfying
Cuntz relations: $s_i^*s_j=\delta _{ij}I, 1\leq i, j\leq n,$ and $\sum s_is_i^*=I.$
It admits many unitarily inequivalent representations. Various
classes of representations of $\CO_n$ have been constructed in
\cite{BJ1, BJ2}, \cite{DKS}. Given a tuple of contractions $\uT=(T_1, \ldots , T_n)$ on
a Hilbert space satisfying $\sum T_iT_i^*=I$, we consider its minimal
isometric dilation $\tilde {\uV }=(\tilde {V_1}, \ldots , \tilde {V_n})$.
We know that the isometries $\tilde {V_i}$ satisfy Cuntz relations and we obtain a
representation $\pi _{\uT}$ of the Cuntz algebra $\CO_n$ by setting
$\pi _{\uT}(s_i)= \tilde {V_i}.$ We wish to classify all representations
of $\CO_n$ we can obtain by dilating {\em commuting\/}  contractive tuples $\uT $.

Let $\CS _n=C(\partial B_n)$ be the $C^*$-algebra of all continuous complex valued
functions on the sphere $\partial B_n=\{ (z_1, \ldots , z_n) : \sum |z_i|^2=1\} .$
We have a distinguished tuple $\underline {z}= (z_1, \ldots , z_n)$
of elements in $\CS _n$
consisting of co-ordinate functions.  Given
any spherical unitary $\uZ= (Z_1, \ldots , Z_n)$
there is a unique representation of  $\CS _n$ which maps
$z_i$ to $Z_i$.  Now given any commuting $n$-tuple of operators
$\uT$, satisfying $\sum T_iT_i^* =I$, we consider its standard commuting
dilation $\tilde {\uS}=(\tilde {S_1}, \ldots , \tilde {S_n})$. Let $\rho _{\uT}$
be the representation of $\CS _n$, obtained by taking $\rho _{\uT}(z_i)=\tilde {S_i}.$

\begin{Definition}{\em
Let $\pi $ be representation of $\CO _n$ on a Hilbert space $\CL $ with
$\uW = (W_1, \ldots , W_n)=(\pi (s_1), \ldots , \pi (s_n))$. The
representation $\pi $ is said to be {\em spherical\/ } if $~\overline
{\mbox {span}}~\{ \uW ^\alpha h: h\in \CL ^c(\uW), \alpha \in \tilde {\Lambda } \}
=\CL $, where $\CL ^c(\uW) $ is the space where the maximal commuting
piece $\uW ^c$ of $\uW$ acts as in Definition 3.}
\end{Definition}

Note that this Definition means in particular that if $\pi $ is spherical
then the maximal commuting piece $\uW ^c$ is non-trivial. We will see
that it is actually a spherical unitary. But this is not a justification
for calling such representations as spherical, because  this happens
for any representation of $\CO _n$, as long as  $\uW ^c$
is non-trivial! The actual justification
of this Definition is in Theorem 18.

\begin{Theorem}
Let $\uT =(T_1, \ldots , T_n)$ be a commuting tuple of operators on a Hilbert
space $\CH $, satisfying $\sum T_iT_i^*=I.$ Then the representation
$\pi _{\uT}$  coming from the minimal isometric dilation of $\uT $ is
spherical. Suppose $\uR  =(R_1, \ldots , R_n)$ is another commuting tuple,
possibly on a different Hilbert space, satisfying $\sum R_iR_i^*=1.$
Then the representations $\pi _{\uT }$, $\pi _{\uR }$ of $\CO _n$ are
unitarily equivalent if and only if the representations $\rho _{\uT }$,
$\rho _{\uR }$ of $\CS _n$ are unitarily equivalent.
\end{Theorem}

\noindent {\sc Proof:} In view of Theorem 13, the maximal commuting piece of the minimal
isometric dilation  $\tilde {\uV }$ of $\uT$ is a
realization of the standard commuting
dilation $\tilde {\uS }$ of $\uT$. The first claim follows easily as
the space on which the standard commuting dilation acts includes the
original space $\CH .$ So $\tilde {\uV }$ is the minimal isometric dilation
of $\tilde {\uS}.$ Similar statement holds for the tuple $\uR$. Now the
Theorem follows due to uniqueness up to equivalence of minimal isometric
dilation of contractive tuples, and unitary equivalence of maximal
commuting pieces of unitarily equivalent tuples.
\qed

So this Theorem reduces the classification problem for representations
of $\CO _n$ arising out of general commuting tuples to that of representations
of $\CS _n$. But $\CS _n$ being a commutative $C^*$-algebra, its representations
are well-understood and is part of standard $C^*$-algebra theory. We find
the description of this theory as presented in Arveson's classic \cite{Ar1}
most suitable for our purposes.

Given any point $w= (w_1, \ldots w_n)\in \partial B_n,$ we have a one dimensional
representation $\phi _{w }$ of  $\CS _n$, which maps $f$ to $f(w)$.
Of course $w $ is a spherical unitary as operator tuple on $\mathbb{C} $.
We can construct the minimal isometric dilation $(W^w_1, \ldots , W^w_n)$
of this tuple as in the proof of
Theorem 15 (Sch\"{a}ffer construction). We see that the dilation space is
$$\CH ^{w}= \mathbb{C}\oplus (\Gamma (\mathbb{C}^n)\otimes \mathbb{C}^n_w)
\subseteq  \mathbb{C}\oplus (\Gamma (\mathbb{C}^n)\otimes \mathbb{C}^n),$$
where $\mathbb{C}^n_w$ is the subspace of vectors orthogonal to $(\overline{w_1},
\ldots ,  \overline{w_n})$ in $\mathbb{C}^n.$ Further the operators
$W^w_i$ are given by
$$W^w_i(h\oplus \sum _{\alpha }e^\alpha \otimes d_{\alpha })
=w_ih\oplus D(e_i\otimes h)\oplus e_i\otimes (\sum _{\alpha } e^\alpha \otimes d_{\alpha }).$$
We denote the associated representation of $\CO _n$
by $\rho _{w }$. This representation is known to be
irreducible as it is  nothing but the GNS representation
of the so-called Cuntz state  on $\CO _n$ (See \cite{DKS}, Example 5.1]), given by
$$s_{i_1}\cdots s_{i_m}s_{j_1}^*\cdots s_{j_p}^*
\mapsto w_{i_1}\cdots w_{i_m}\overline {w_{j}}_1\cdots \overline {w_{j}}_p.$$

Now an arbitrary multiplicity free representation of $\CS_n$ can be
described as follows \cite{Ar1}. Consider a finite Borel measure $\mu $ on
$\partial B_n$. Then we get a representation of $\CS _n$ on the
Hilbert space $L^2(\partial B_n, \mu)$, which sends $f\in \CS _n$ to
the operator `multiplication by $f$'. This representation can be
thought of as direct integral of representations $\phi _{w}$
with respect to measure $\mu $. Now it is not hard to see that the
associated representation of $\CO _n$ is simply the direct
integral of representations $\rho _{w}$ with respect to measure
$\mu $ and acts on $\int \! \! \!\! \!  \oplus\CH ^{w}\mu (dw).$
Finally an arbitrary representation of $\CS _n$ is a countable
direct sum of such multiplicity free representations. So we have
proved the following result.

\begin{Theorem}
Every spherical representation of $\CO _n$ is a direct integral
of representations $\rho _w, w\in \partial B_n$ (GNS representations
of Cuntz states).
\end{Theorem}

Here we have not bothered to write down as to when two such representations
are equivalent. But in view of Theorem 17, we can do it exactly as in
 (\cite{Ar1}, page:54-55), by keeping track of multiplicities and equivalence
classes of measures.

\begin{Theorem}
Let $\pi $ be a representation of $\CO _n$. Then
(i) $\pi $ decomposes uniquely as $\pi=\pi ^0\oplus \pi ^1$,
where $\pi ^0$ is spherical and $(\pi ^1 (s_1), \ldots , \pi ^1 (s_n))$
has trivial maximal commuting piece (Either $\pi ^0$ or $\pi ^1$ could
also be absent); (ii) The maximal commuting piece
of $(\pi (s_1), \ldots , \pi (s_n))$ is either trivial or it is a
spherical unitary. (iii) If $\pi $ is irreducible then either the maximal
commuting piece is trivial or it is one dimensional. In the second case,
it is unitarily equivalent to GNS representation of a Cuntz state.
\end{Theorem}

\noindent {\sc Proof:} Suppose $\pi $ is a representation of $\CO _n$ on a Hilbert
space $\CL $ and $\uW =(\pi (s_1), \ldots , \pi (s_n)).$ Consider the
space $\CL ^0$ generated by $\CL ^c(\uW)$ as
$\CL ^0=~\overline{\mbox{span}}~\{ \uW ^\alpha h: h\in \CL ^c(\uW), \alpha
\in \tilde{\Lambda}\}.$ Now each $W_i^*$ leaves $\CL ^c(\uW)$ invariant
and clearly $\CO _n=C^*\{\us ^\alpha (\us ^\beta)^*: \alpha , \beta \in \tilde
{\Lambda }\}.$ Then it follows that $\CL ^0$ is a reducing subspace for $\pi .$
Taking $\CL ^1=(\CL ^0)^\perp $, we decompose $\pi $ as $\pi ^0\oplus
\pi ^1$ with respect to $\CL =\CL ^0\oplus \CL ^1.$ It is clear that this
is a decomposition as required by (i). Uniqueness of this decomposition
and (ii) follow easily as maximal commuting piece of direct sum of tuples
is direct sum of maximal commuting pieces (Corollary 5) and then (iii)
follows from Theorem 18.
\qed

Let us see as to what happens if we dilate commuting tuples $\uT $, satisfying
just $\sum T_iT_i^*\leq I$. In this case, as is well-known, the minimal
isometric dilation decomposes as $((V_1\otimes I)\oplus W_1,
\ldots , (V_n\otimes I)\oplus W_n)$ where $(V_1, \ldots , V_n)$ is the
standard tuple of full Fock space, and $(W_1, \ldots , W_n)$ are isometries
satisfying Cuntz relations. If $\uT $ is not pure the term $(W_1, \ldots,
W_n)$ is present and we get a representation of  $\CO _n$. However, as
seen in the proof of Theorem 13, $(W_1, \ldots , W_n)$ is a minimal
isometric dilation of a spherical tuple $(Z_1, \ldots , Z_n)$ (the
`spherical part' of the standard commuting dilation of $\uT$) and hence
the representation of $\CO _n$ we get is still spherical.

Finally we remark that it is easy to get examples of
non-commuting tuples dilating to representations of $\CO _n$ which
are not spherical. For instance we
can consider  the tuple $\uR=(R_1, R_2)$ on ${\mathbb {C}}^2$ defined by
$$R_1=\left[ \begin{array}{cc}0& 1\\0&0\end{array}\right],
R_2=\left[ \begin{array}{cc}0& 0\\1&0\end{array}\right].$$
Then as $R_1R_1^*+R_2R_2^*=1$, the minimal isometric dilation of
$(R_1, R_2)$ satisfies Cuntz relations. We can see that it has
trivial commuting piece through a simple application of  Corollary 4.3 of
\cite{DKS}.

\vsp
\noindent{\bf Acknowledgements:}
 The first author is supported by Indo-French (IFCPAR) grant No: IFC/ 2301-1/2001/1253,
the second author is supported by research grant no. SR/FTP/MS-16/2001 of the Department of
Science and Technology, India and the last author is supported by
a research fellowship from the Indian Statistical Institute.

\end{section}

\begin{tabular}{ll}
 {\sc B. V. Rajarama Bhat and Santanu Dey}& ~~~~~ ~~~~~ {\sc Tirthankar Bhattacharyya}\\
  Indian Statistical Institute,& ~~~~~ ~~~~~ Department of Mathematics,\\
R. V. College Post,  & ~~~~~ ~~~~~  Indian Institute of Science \\
Bangalore 560059, India. & ~~~~~ ~~~~~ Bangalore 560012, India.\\
 e-mail: {\sl bhat@isibang.ac.in} & ~~~~~ ~~~~~  {\sl tirtha@math.iisc.ernet.in}\\
and {\sl santanu@isibang.ac.in }\\
\end{tabular}

\end{document}